\documentclass[11pt]{amsart}
\usepackage{graphicx}
\usepackage{amsmath,amsfonts,amssymb,latexsym,amsthm,amscd}
 \newtheorem{theorem}{Theorem}[section]
 \newtheorem{lemma}[theorem]{Lemma}
 \newtheorem{proposition}[theorem]{Proposition}

\numberwithin{equation}{section}
\begin{document}

\title[self mapping degrees of spherical manifolds]
  {On Self-mapping Degrees of $S^3$-geometry manifolds}

\author{Xiaoming Du}
\address{School of Mathematical Sciences, Peking University,
  Beijing 100871, P.R.China}
\email{duxiaoming@math.pku.edu.cn}

\begin{abstract}
In this paper we determined all of the possible self mapping degrees
of the manifolds with $S^3$-geometry, which are supposed to be all
3-manifolds with finite fundamental groups. This is a part of a
project to determine all possible self mapping degrees of all
closed orientable 3-manifold in Thurston's picture.\\[4mm]
\end{abstract}

\thanks{I would like to thank Professor
Wang, Shicheng for drawing my attention to this topic. The author is
partially supported by grant No.10631060 of the National Natural
Science Foundation of China and Ph.D. grant No. 5171042-055 of the
Ministry of Education of China.}

\maketitle

\section{Introduction and Results}

For a closed oriented manifold $M$, deciding the mapping-degree set
$D(M)=\{\,\text{deg}f:\;|\,f:M\to\,M\}$ is a natural problem. When
the dimension is 1 and 2, the answer is well-known. When the
dimension is higher than 3, there are many interesting results in
this topic. For classical discussions see \cite{Olum1}, \cite{Olum2}
and for quite recent papers see \cite{Duan}, \cite{Ding} and
\cite{Tan}. But it is difficult to get general results, since there
is no classification result for manifolds of dimension $n>3$.

The case of dimension 3 becomes attractive in the topic and it
becomes possible to calculate $D(M)$ for any closed oriented
3-manifold $M$. Since Thurston's geometrization conjecture, which
seems to be confirmed, implies that closed oriented 3-manifolds can
be classified in a reasonable sense.

Thurston's geometrization conjecture claims that the each
Jaco-Shalen-Johanson decomposition piece of a prime 3-manifold
supports one of eight geometries, which are $H^3$, $\widetilde
{PSL}(2,\mathbb{R})$, $H^2\times E^1$, $Sol$, $Nil$, $E^3$, $S^3$
and $S^2\times S^1$ (\,for details see \cite{Thurston} and
\cite{Scott}). Call a closed orientable 3-manifold $M$ is {\it
geometric} if it supports one of the eight geometries above.

Hence we should ask first

\textbf{Question(*)}\;\,How to determine $D(M)$ for a geometric
3-manifold $M$?

Indeed, for most cases the answer to this question is known. When
$M$ supports the geometry of $H^3$ or
$\widetilde{PSL}(2,\mathbb{R})$, by using the Gromov volume or
$\widetilde{PSL}(2,\mathbb{R})$ volume, $D(M)$ should be either
$\{\,0,1,-1\}$ or $\{\,0,1\}$, depending on whether $M$ admits a
self map of degree $-1$ or not (\,see \cite{Wang1} for details). For
3-manifolds admitting $E^3$ and $Sol$ structures and for some $Nil$
3-manifolds, $D(M)$'s have been determined just recently in
\cite{Sun}.

The remaining important case for question(*) is for the 3-manifolds
supporting $S^3$-geometry, which are supposed to be all 3-manifolds
with finite fundamental groups.

The question(*) for this case has several partial answers before.
The answer was known for Poincare homology 3-spheres in
\cite{Plotnik}, and for Quaternion spaces in \cite {Hayat2}. More
generally, for the degrees of the self-mappings inducing
automorphisms of $\pi_1(M)$, the answer had been known in
\cite{Hayat4}, which stated that $\{\text{deg}f\,|\;
f_*:\pi_1(M)\to\pi_1(M)\;\text{is an automorphism}\}=
\{\,k^2\,|\,(\,k,|\pi_1(M)|)=1\}+|\pi_1(M)|\cdot\mathbb{Z}$\;
(\cite{Hayat4}, Theorem 2.2).

In this paper, we will determine the degrees of the mappings
inducing all possible endomorphisms of the fundamental groups to
give a complete answer of question(*) for every 3-manifold $M$
admitting $S^3$-geometry.

\bigskip

According to \cite{Scott}, the fundamental group of a 3-manifold
with $S^3$-geometry structure is among the following eight types:
$\mathbb{Z}_p$\,, $D_{4n}^*$\,, $T_{24}^*$\,, $O_{48}^*$\,,
$I_{120}^*$\,,
$T_{8\cdot\,3^{^q}}^{'}$\,,\,$D_{n'\cdot\,2^{^q}}'\;(\,2\nmid n'\,)$
and $\mathbb{Z}_m\times G$ where $G$ belongs to the previous seven
ones and $|G|$ is coprime to $m$. $T_{24}^*$\,, $O_{48}^*$ and
$I_{120}^*$ have intuitive explanations: they are the pre-images of
the rigid rotation groups of the regular tetrahedron, octahedron and
icosahedron under the double covering $S^3\to SO(3)$. $D_{4n}^*$ is
the pre-image of the dihedral group $D_{2n}$ under the same double
covering. $T_{8\cdot\,3^{^q}}^{'}$\,,\,$D_{n'\cdot\,2^{^q}}'$ and
$\mathbb{Z}_m\times G$ are discrete subgroups of $SO(4)$ which act
on $S^3$. These groups have presentations as following:

\begin{table}[h]
\begin{center}
\renewcommand{\arraystretch}{1.3}
\begin{tabular}{|c|c|}
  \hline
  $\pi_1(M)$ & presentation\\
  \hline
  $\mathbb{Z}_p$ & $\langle\,a\,|\,a^p=1\,\rangle$\\
  \hline
  $D_{4n}^*$ &
  $\langle\,a,b\;|\;a^2=b^n=(ab)^2, a^4=1\rangle$\\
  \hline
  $T_{24}^*$ & $\langle\,a,b\,|\,a^2=b^3=(ab)^3, a^4=1\rangle$\\
  \hline
  $O_{48}^*$ & $\langle\,a,b\,|\,a^2=b^3=(ab)^4, a^4=1\rangle$\\
  \hline
  $I_{120}^*$ & $\langle\,a,b\,|\,a^2=b^3=(ab)^5, a^4=1\rangle$\\
  \hline
  $T_{8\cdot3^q}^{'}$ &
    $\begin{array}{ll} & \langle\,i,j,k,w\,|\,i^2=j^2=k^2, i^4=1, ij=k, jk=i, ki=j, w^{3^q}=1,\\
      &\qquad w\cdot\,i\cdot\,w^{-1}=j\,;\;\;w\cdot\,j\cdot\,w^{-1}=k\,;\;\;w\cdot\,k\cdot\,w^{-1}=i\,\rangle
    \end{array}$\\
  \hline
  $D_{n'\cdot\,2^q}'$ &
    $\langle\,w,u\,|\,u^{n'}=1, w^{2^q}=1, w\cdot\,u\cdot\,w^{-1}=u^{-1}\rangle$\\
  \hline
  $\mathbb{Z}_m\times G$ &
    $G\;\,\text{is in one of the seven types above and}\;(m,|G|)=1$\\
  \hline
\end{tabular}
\end{center}
\end{table}

Respectively, we have the following theorem on the self-mapping
degrees.

\begin{theorem}
The sets of self-mapping degrees of $S^3$-geometry manifolds are
listed as follow:
\end{theorem}

\begin{table}[h]
\begin{center}
\renewcommand{\arraystretch}{1.3}
\begin{tabular}{|c|c|}
  \hline
  $\pi_1(M)$ & $D(M)$\\
  \hline
  $\mathbb{Z}_p$ & $\{\,k^2\,|\,k\in\mathbb{Z}\}+p\,\mathbb{Z}$\\
  \hline
  $D_{4n}^*$ &
  $\{\,h^2\,|\,h\in\mathbb{Z};\,2\nmid\,h\;\text{or}\;h=n\,\text{or}\;h=0\}+4n\mathbb{Z}$\\
  \hline
  $T_{24}^*$ & $\{\,0,1,16\}+24\mathbb{Z}$\\
  \hline
  $O_{48}^*$ & $\{\,0,1,25\}+48\mathbb{Z}$\\
  \hline
  $I_{120}^*$ & $\{\,0,1,49\}+120\mathbb{Z}$\\
  \hline
  $T_{8\cdot3^q}^{'}$ & $\left\{\begin{array}{ll}
      \{\,k^2\cdot\,(3^{2q-2p}-3^{q})\;\;|\,3\nmid\,k, q\geq\,p>0\}+8\cdot3^q\mathbb{Z}&(\,2\,|\,q\,)\\
      \{\,k^2\cdot\,(3^{2q-2p}-3^{q+1})\;\;|\,3\nmid\,k,q\geq\,p>0\}+8\cdot3^q\mathbb{Z}&(\,2\nmid\,q\,)
      \end{array}\right.$\\
  \hline
  $D_{n'\cdot\,2^q}'$ &
    $\begin{array}{ll}&\{\,k^2\cdot[\,1-(n')^{2^{^{\,q}}-1}\,]^i\cdot
      [\,1-2^{\,(\,2p-q\,)(\,n'-1\,)}\,])^j\;|\,i,j,k,p\in\mathbb{Z},\\
      &\qquad\qquad\qquad\,q\geq\,p>0\}+n'\,2^{\,q}\,\mathbb{Z}\end{array}$\\
  \hline
  $\mathbb{Z}_m\times G$ & $\left\{d\in\mathbb{Z}\,\left| \begin{array}{l}
    d\equiv\,h\:(\text{mod\;}|G|)\;(h\in\,D(N),N\:\text{is the}\\
      \qquad S^3\text{-geometry manifold with }\pi_1(N)=G)\\
    d\equiv\,k^2
    (k\in\mathbb{Z})(\text{mod}\;m)\end{array}\right.\right\}$\\
  \hline
\end{tabular}
\end{center}
\end{table}

The proof of this theorem will be divided into two parts. Firstly,
since Proposition 2.1\,(see bellow) states that for each spherical
3-manifold, there is a well defined mapping from the set of
endomorphisms of the fundamental group to the mod $|\pi_1|$ classes
of the self-mapping degrees, we need to find out all possible
endomorphisms of the fundamental groups of the above eight types.
Secondly, for each endomorphism, we need to calculate the degree of
some self-mapping realizing it.

\section{Preliminaries}

Since the second homotopy group of a $S^3$-geometry manifold is
trivial, the existences of self-mappings can be detected by the
obstruction theory. P. Olum showed in \cite{Olum1} the first and in
\cite{Olum2} the second part of the following proposition.

\begin{proposition}[Olum]
Let $M$ be an orientable 3-manifold with finite fundamental group
and trivial $\pi_2(M)$. Every endomorphism $\phi:\pi_1(M)\to
\pi_1(M)$ is induced by a (basepoint preserving) continuous map
$f:M\to M$. Furthermore, if $g$ is also a continuous self-mapping of
$M$ such that $f_{*}=g_{*}=\phi$ then $\deg{f}\equiv\deg{g}\mod |\pi_1(M)|$.
\end{proposition}

According to this proposition, the self-mapping degrees of a
spherical 3-manifold $M$ are closely related to the endomorphisms of
$\pi_1(M)$.

The statement and proof of the next lemma is elementary (see
\,\cite{Hayat4} lemma 3.4). This lemma is very useful in the
calculation of the degree of a self-mapping corresponding to a given
endomorphism on the fundamental group.

\begin{lemma}[Hayat-Kudryavtseva-Wang-Zieschang]
  Let $f$ be a self-mapping of $M$, $G=\pi_1(M)$ contain a
  subgroup $H$ such that $f_*(H)\subset gHg^{-1}$ for some $g\in G$.
  Consider the covering $p:\widetilde{_HM}\to M$ corresponding to $H$, that is,
  $H=p_*(\pi_1(\widetilde{_HM}))$. Then there is a map
  $\widetilde{_Hf}:\widetilde{_HM}\to\widetilde{_HM}$ and a
  homeomorphism $J:M\to M$ isotopic to $id_M$ such that the
  following diagram is commutative:
  $$\CD
  \widetilde{_HM} @> \widetilde{_Hf} >> \widetilde{_HM} \\
  @V p\; VV @V VV J\circ p  \\
  M @>f>> M.
  \endCD$$
  A consequence is that \text{deg}($\widetilde{_Hf}$)=\text{deg} $f$.
\end{lemma}

If $H$ is a Sylow subgroup of $\pi_1(M)$, $H$ satisfies the
condition of the above lemma by the second Sylow theorem
(\cite{Jacobson}, 1.13).

\section{Proof of Theorem 1.1}

At first we will find out all the possible endomorphisms of the
fundamental groups of the spherical 3-manifolds.
\begin{lemma}
Two endomorphisms of the fundamental group of a spherical 3-manifold
having isomorphic kernels can be transformed to each other by a
composition with an isomorphism of the fundamental group. The
kernels of the non-trivial endomorphisms of the fundamental groups
of the spherical 3-manifolds are listed as follow:
\end{lemma}

\begin{table}[h]
\begin{center}
\renewcommand{\arraystretch}{1.3}
\begin{tabular}{|c|c|}
  \hline
  $\pi_1(M)$ & non-trivial kernels of endomorphisms of $\pi_1(M)$\\
  \hline
  $\mathbb{Z}_p$ & $\mathbb{Z}_{p/k}$\\
  \hline
  $D_{4n}^*$ & $D_{4\frac{n}{2}}^*$\;(if $2\,|\,n$), $\mathbb{Z}_{2n}$ and $\mathbb{Z}_h$\;($2 \nmid h$)\\
  \hline
  $T_{24}^*$ & $Q_8$\\
  \hline
  $O_{48}^*$ & $T_{24}^*$\\
  \hline
  $I_{120}^*$ & none\\
  \hline
  $T_{8\cdot3^q}^{'}=\,Q_8\rtimes\mathbb{Z}_{3^q}$ & $Q_8\times\mathbb{Z}_{3^{^{q-p}}}$\;(\,$q\geq\,p>0$)\\
  \hline
  $D_{n'\cdot\,2^q}'=\,\mathbb{Z}_{n'}\rtimes\mathbb{Z}_{2^q}$ &
    $\mathbb{Z}_{\,n'}\times\mathbb{Z}_{\,2^{^{\,q-p}}}$\;(\,$q\geq\,p>0$\,) and
    $\mathbb{Z}_{n'/n''}\times\mathbb{Z}_{\,2^{^{\,q}}}\;(\,n''\,|\,n', n''>1)$\\
  \hline
  $\mathbb{Z}_m\times G$ & $\mathbb{Z}_{m'}\times H$\;($m'\,|\,m, H\subset\,G$)\\
  \hline
\end{tabular}
\end{center}
\end{table}

\begin{proof}

\textbf{Case I:}

Since $\mathbb{Z}_p=\langle\,a\,|\,a^p=1\,\rangle$ is a cyclic
group, the result is obvious. The endomorphism is like
$\phi:a\mapsto\,a^t$ and its kernel is
$\langle\,a^{\frac{p}{(t,p)}}\rangle\cong\mathbb{Z}_{(t,p)}$.

\textbf{Case II:}

$D_{4n}^*$ has a representation as
$\langle\,a,b\;|\;a^2=b^n=(ab)^2=-1\rangle$,~$n=2^qn'$,~$n'$~is an
odd number (see \cite{Hayat4}).

First note that $ab^ka^{-1}=b^{2n-\,k}$, so
$\langle\,b^k\rangle\cong\mathbb{Z}_{\frac{2n}{(2n,\:k)}}$.
Secondly, The elements in $D_{4n}^*\setminus\langle\,b\,\rangle$ are
like $ab^l$ and have order 4. We will discuss the normal subgroup
according to whether it includes such elements or not.

\noindent (I)~If the normal subgroup $H$ includes some element in
$D_{4n}^*\setminus\langle\,b\,\rangle$, because
$b(ab^l)b^{-1}=ab^{l-2}$, we have $b^2\in\,H$,
$H\supseteq\langle\,ab^l,b^2\rangle$. In fact,
$\langle\,ab^l,b^2\rangle$ is a proper subgroup of $D_{4n}^*$ if and
only if $2\,|\,n$. When $2\,|\,n$, we have
$\langle\,ab^l,b^2\rangle\,\cong\,D_{4\cdot\frac{n}{2}}^*$ and
$D_{4n}^*/\langle\,ab^l,b^2\rangle$~$=D_{4n}^*/D_{4\frac{n}{2}}^*$~$\cong\mathbb{Z}_2\leq\,D_{4n}^*$.
The corresponding endomorphism is like
$\phi:a\mapsto(-1)^l,\;b\mapsto-1$.

\noindent (II)~If the normal subgroup~$H$~is only generated by some
$b^k$, then

(A)~When $2\,|\,\frac{2n}{(2n,\:k)}$, we have
$-1\in\langle\,b^k\rangle$. The quotient group
$D_{4n}^*/\langle\,b^k\rangle$ is homomorphic to the dihedral group
$D_{2(2n,\:k)}$ and also a subgroup of $D_{4n}^*$. However, there
are ~$(2n,\:k)$ elements of order 2 in~$D_{2(2n,\:k)}$ while there
is only one in $D_{4n}^*$. This is absurd when $(2n,\,k)>1$. So the
kernels of endomorphisms of $D_{4n}^*$ cannot be
~$\langle\,b^k\rangle$ with $(2n,\,k)>1$. When $(2n,\,k)=1$,
$\langle\,b^k\,\rangle=\langle\,b\,\rangle$, we have
$D_{4n}^*/\langle\,b\,\rangle=D_{4n}^*/\mathbb{Z}_{2n}\cong\mathbb{Z}_2\leq\,D_{4n}^*$.
So there exists only one endomorphism of $D_{4n}^*$ the kernel of
which is $\langle\,b\,\rangle\cong\mathbb{Z}_{2n}$. This
endomorphism is like $\phi:a\mapsto-1,\;b\mapsto1$.

(B)~When $2\nmid\frac{2n}{(2n,\:k)}$~(i.e.~$2^{q+1}\,|\,k$), we have
$-1\not\in\langle\,b^k\rangle$ and $D_{4n}^*/\langle\,b^k\rangle$~
$=D_{4n}^*/\mathbb{Z}_{\frac{2n}{(2n,\:k)}}$
$\cong$~$D_{4\frac{(2n,\,k)}{2}}^*\leq\,D_{4n}^*$. So there exists
an endomorphism of $D_{4n}^*$ the kernel of which is
$\langle\,b^k\rangle$. The endomorphisms having this kernel are like
$\phi:a\mapsto\,b^sab^{-s},\;b\mapsto\,b^{\frac{2n}{(2n,\:k)}\cdot\,t}\,
(\,s,t \in \mathbb{Z},(t,k)=1)$, which are conjugate to
$\psi:a\mapsto\,a,\;b\mapsto\,b^{\frac{2n}{(2n,\:k)}\cdot\,t}$. By
the composition with the isomorphism $a \mapsto a,\;b \mapsto b^t$
of $D_{4n}^*$, they can be transformed to
$\phi_0:a\mapsto\,a,\;b\mapsto\,b^{\frac{2n}{(2n,\:k)}}$.

\textbf{Case III:}

$T_{24}^*$ has a representation as
$\langle\,a,b\,|\,a^2=b^3=(ab)^3=-1\rangle$\;.

The normal subgroups of $T_{24}^*$ are:\\
(I)\;$\{\,1,-1\}$.\; But $T_{24}^*/\{\,1,-1\}\cong\,T_{12}=A_4$ is not a subgroup of $T_{24}^*$;\\
(II)\;$\langle\,a,bab^{-1}\rangle\cong\,Q_8$. In this case
$T_{24}^*/\langle\,a,bab^{-1}\rangle\cong\mathbb{Z}_3$ and subgroups
of $T_{24}^*$ isomorphic to $\mathbb{Z}_3$ have the forms like
$\{1,xb^2x^{-1},xb^4x^{-1}\}$.

$T_{24}^*$ has no other normal subgroup. In fact, if $H$ is a normal
subgroups of $T_{24}^*$ including $a^2=-1$, then $H/\{\,1,-1\}$ is
isomorphic to a normal subgroup of the tetrahedron group
$T_{12}\cong\,A_4$. But the only normal subgroup of $T_{12}$ is
$\mathbb{Z}_2\oplus\mathbb{Z}_2$. Its pre-image in $T_{24}^*$ is
$Q_8$. If $H$ does not include $-1$, then $H$ must include an
element of $T_{24}^*$ of odd order. The odd-order elements of
$T_{24}^*$ conjugate to $b^2$. So $b^2\in\,H,\,ab^2a^{-1}b^2\in\,H$.
But $ab^2a^{-1}b^2$ is an element of order six.~It is conjugated to
$b$. This implies~$b\in\,H,\,-1=b^3\in\,H$, contradiction.

The endomorphisms corresponding to the normal subgroup $Q_8$ are
like $\phi:a\mapsto1,\;b\mapsto\,xb^2x^{-1}$. They are conjugate to
$\phi_0:a\mapsto1,\;b\mapsto\,b^2$.

\textbf{Case IV:}

$O_{48}^*$ has a representation as
$\langle\,a,b\,|\,a^4=b^3=(ab)^2=-1\rangle$\;.

The normal subgroups of $O_{48}^*$ are:\\
(I)\;$\{1,-1\}$.\; But $O_{48}^*/\{\,1,-1\}\cong\,O_{24}=S_4$ is not a subgroup of $O_{48}^*$;\\
(II)\;$\langle\,a^2,b\,\rangle\cong\,T_{24}^*$. In this case
$O_{48}^*/T_{24}^*\cong\mathbb{Z}_2\leq\,O_{48}^*$.

Similar to the last case, $O_{48}^*$ has no other normal subgroup.

The endomorphism corresponding to the normal subgroup $T_{24}^*$ is
like $\phi:a\mapsto-1,\;b\mapsto1$.

\textbf{Case V:}

$I_{120}^*=\langle\,a,b\,|\,a^2=b^3=(ab)^5=-1\rangle$\;.

The only non-trivial normal subgroup of $I_{120}^*$ is
$\{\,1,\,-1\}$. But $I_{120}^*/\{\,1,\,-1\}\cong\,I_{60}=A_5$ which
is not a subgroup of $I_{120}^*$. So the $I_{120}^*$ only admits
trivial endomorphism.

\textbf{Case VI:}

$T_{8\cdot3^q}^{'}=\,Q_8\rtimes\mathbb{Z}_{3^q}$,\,
$Q_8=\{\,\pm\,1,\,\pm\,i,\,\pm\,j,\,\pm\,k\,\},\,\mathbb{Z}_{3^q}=\langle\,w\,\rangle$,
$w\cdot\,i\cdot\,w^{-1}=j\,;\;\;w\cdot\,j\cdot\,w^{-1}=k\,;\;\;w\cdot\,k\cdot\,w^{-1}=i.$

If $K\unlhd\,T_{8\cdot\,3^q}'$, for the symmetry of $i$, $j$ and
$k$, suppose $w^{\alpha}i\in\,K$, then
$w^{-1}(w^{\alpha}i)w=w^{\alpha-1}(wk)=w^{\alpha}k\in\,K$. So
$i^{-1}k=j\in\,K$, similarly $k,i\in\,K$ and $Q_8\subset\,K$.

Since $w^3$ commutes with the elements in $Q_8$, the normal
subgroups of $T_{8\cdot\,3^q}'$ should be like
$Q_8\times\langle\,w^{3^{^{p}}}\,\rangle\cong
Q_8\times\mathbb{Z}_{3^{q-p}}\;(\,q\geq\,p>0\,)$ and
$(Q_8\rtimes\mathbb{Z}_{3^{\,q}})/(Q_8\times\mathbb{Z}_{3^{q-p}})$
$\cong\mathbb{Z}_{\,3^{^{p}}}$
$\cong\langle\,w^{{3^{^{q-p}}}}\,\rangle$ which is isomorphic to a
subgroup of $T_{8\cdot3^q}^{'}$.

The endomorphisms corresponding to the normal subgroup
$Q_8\times\langle\,w^{3^{^{p}}}\,\rangle\cong
Q_8\times\mathbb{Z}_{3^{^{q-p}}}$ are like
$\phi:i\mapsto1,\,j\mapsto1,\,k\mapsto1,w\mapsto\,w^{h\cdot\,3^{^{q-p}}}$.

\textbf{Case VII:}

$D_{n'\cdot\,2^q}'=\,\mathbb{Z}_{n'}\rtimes\mathbb{Z}_{2^q}$,
$2\nmid n'$,
$\mathbb{Z}_{n'}=\langle\,u\,\rangle,\,\mathbb{Z}_{2^q}=\langle\,w\,\rangle,\;w\cdot\,u\cdot\,w^{-1}=u^{-1}.$

Any element in $\mathbb{Z}_{n'}\rtimes\mathbb{Z}_{2^q}$ can be
written as $w^{\alpha}u^{\beta}$. If $K$ is a normal subgroup, for
$w^{\alpha}u^{\beta}\in\,K$, we have
$w^{-1}(w^{\alpha}u^{\beta})w$~$=w^{\alpha}(w^{-1}u^{\beta}w)$~$=w^{\alpha}u^{-\beta}\in\,K$
and
$u^{\beta}w^{\alpha}=u^{\beta}(w^{\alpha}u^{\beta})u^{-\beta}\in\,K$.
So
$u^{2\beta}\in\,K$,~$u^{(\,2\beta,\,n')}=u^{(\,\beta,\,n')}\in\,K$,
and then $u^{\beta}\in\,K$,~$w^{\alpha}\in\,K$.

If some $w^{\alpha}u^{\beta}$ satisfies $2\nmid\alpha$, then
$w=w^{(\alpha,\,2^{^q})}\in\,K$,
$u(w^{\alpha}u^{\beta})u^{-1}$~$=w^{\alpha}u^{(-1)^{^{\alpha}}}u^{\beta-1}$~$=w^{\alpha}u^{\beta-2}\in\,K$.
Hence $u^2\in\,K$ and $u=u^{(\,2,\,n')}\in\,K$. Then
$K=\langle\,u,w\,\rangle$ is a trivial subgroup.

If every element $w^{\alpha}u^{\beta}$ in $K$ satisfies
$2\,|\,\alpha$, let $\alpha_{_0}$ be the greatest common divisor of
such $\alpha$'s. We have $w^{\,\alpha_{_0}}\in\,K$,
$\langle\,w^{(\,\alpha_{_0},\,2^{^q})}\,\rangle=\langle\,w^{2^{^p}}\rangle\subset\,K\;(\,q\geq\,p>0\,)$.
Note that $w^{2^{^p}}$ commutes with the elements in
$\mathbb{Z}_{n'}\rtimes\,\mathbb{Z}_{2^{\,q}}$. Hence
$K\cong\langle\,u^{n''}\,\rangle\times\langle\,w^{2^{^p}}\,\rangle$
$\cong\mathbb{Z}_{\,n'/n''}\times\mathbb{Z}_{\,2^{^{\,q-p}}}\;(n''|n')$
and
$(\mathbb{Z}_{n'}\rtimes\mathbb{Z}_{2^{\,q}})/K\cong\mathbb{Z}_{n''}\rtimes\mathbb{Z}_{2^p}$.
This quotient group is isomorphic to a subgroup of
$\mathbb{Z}_{n'}\rtimes\mathbb{Z}_{2^{\,q}}$ if and only if
(A)\,$n''=1$ and $q\geq\,p>0$ or (B)\,$n''>1$ and $p=0$. So
endomorphisms $\phi$ of $\mathbb{Z}_{n'}\rtimes\mathbb{Z}_{2^{\,q}}$
have the forms
(A)\,$\phi(u)=1,\,\phi(w)=w^{h\cdot\,2^{^{\,q-p}}}\;(\,2 \nmid h,
q\geq\,p>0)$ or
(B)\,$\phi(u)=u^{l\cdot\,n'/n''},\,\phi(w)=1\;(n''>1, (n'',l)=1)$.

\textbf{Case VIII:}

Since $\mathbb{Z}_m\times\,G$ is a direct product, the result is
obvious.
\end{proof}

In the following part we will calculate the degree of the
self-mappings corresponding to the above endomorphisms. For the
cases $\pi_1(M)=\mathbb{Z}_p, D_{4n}^*$ and $T_{24}^*$, by using the
coordinate system of $SU(2)$, we construct the concrete mappings.
For the remaining cases, by using lemma 2.2, we reduce the problem
to the known cases.

We denote
$S^3=SU(2)=\{(z_1,z_2)\in\mathbb{C}^2\;|\;|z_1|^2+|z_2|^2=1\}$.

\textbf{Case I:}

When $M$ is a lens space and
$\pi_1(M)=\mathbb{Z}_p=\langle\,a\,|\,a^p=1\,\rangle$, the action of
$\mathbb{Z}_p$ on $S^3$ can be realized as
$a\circ(z_1,z_2)=(\zeta\,z_1,\zeta^q\,z_2)$, here $\zeta$ is an
$p^{th}$ root of 1.

Let $\phi:\mathbb{Z}_p\to\mathbb{Z}_p\,, a\mapsto a^k$ be an
endomorphism with the arbitrary integer $k$. Then take
${\widetilde{f}}':S^3\to\,\mathbb{C}^2,(z_1,z_2)\mapsto
(z_1^k,z_2^k)$ and $\widetilde{f}=\pi\circ \widetilde{f}':S^3\to
S^3$, here $\pi:\mathbb{C}^2\to S^3$ is the standard radial
projection.

It is easy to check $\widetilde{f}\circ a=a^k\circ\widetilde{f}$.
This means we can well-define $f:M\to M$ with
$f_*=\phi:\pi_1(M)\to\pi_1(M)$ and
~deg~$f=$~deg~$\widetilde{f}=k^2$.

So
$$D(L_{p,\,q})=\{\,k^2\,|\,k\in\mathbb{Z}\}+p\,\mathbb{Z}.$$

\textbf{Case II:}

$\pi_1(M)=D_{4n}^*=\langle\,a,b\;|\;a^2=b^n=(ab)^2=-1\rangle$,~$n=2^qn'$,~$n'$~is
an odd number.

$D_{4n}^*$ is the double cover of the dihedral group $D_{2n}\subset
SO(3)$ in $S^3$. To calculate the endomorphisms and construct the
corresponding mappings, take a representation
$\rho:D_{4n}^*\to\,SU(2)$ as

$a\,\mapsto\left(
\begin{array}{ccc}
   0 & 1\\
  -1 & 0
\end{array} \right),
b\,\mapsto\left( \begin{array}{ccc}
  \zeta & 0\\
      0 & \bar{\zeta}
\end{array} \right)$,~here~$\zeta^n=-1$.

Then the left actions of elements of $D_{4n}^*$ on $SU(2)$ are
$a\,\circ\,(z_1,z_2)=(-\bar{z}_2,\bar{z}_1)$,
$b\,\circ\,(z_1,z_2)=(\zeta\,z_1,\zeta\,z_2)$.

When the self-mapping induces a endomorphism having the kernel as
$D_{4\cdot\frac{n}{2}}^*\cong\langle\,a,b^2\,\rangle\;(2\,|\,n)$,
take
$${\widetilde{f}}':S^3\to\mathbb{C}^2,(z_1,z_2)\mapsto\,
(z_1^n+\bar{z}_2^{\,n},z_1z_2^{\,n-1}-\bar{z}_1^{\,n-1}\bar{z}_2)$$
and let $\widetilde{f}=\pi\,\circ{\widetilde{f}}':S^3\to\,S^3$, here
$\pi:\mathbb{C}^2\to\,S^3$ is the standard radial projection. Then
$\widetilde{f}$ satisfy
$\widetilde{f}\circ\,a=\widetilde{f},\widetilde{f}\circ\,b=b^n\circ\widetilde{f}$.
This induces $f:S^3/D_{4n}^*\to\,S^3/D_{4n}^*$ with
$kerf_*\cong\,D_{4\cdot\frac{n}{2}}^*$ and deg$f$=~$n^2$.

When the self-mapping induces a endomorphism having the kernel as
$\mathbb{Z}_{2n}$, take
$${\widetilde{f}}':S^3\to\,\mathbb{C}^2,(z_1,z_2)\mapsto
(z_1^{2n}-\bar{z}_2^{2n},\,z_1z_2^{2n-1}+\bar{z}_1^{2n-1}\bar{z_2})$$
and let $\widetilde{f}=\pi\,\circ{\widetilde{f}}':S^3\to\,S^3$. Then
$\widetilde{f}\circ\,a$ $=b^n\circ\,\widetilde{f}$
$=-\widetilde{f},\widetilde{f}\circ\,b$ $=\widetilde{f}$, which
induces $f:S^3/D_{4n}^*\to\,S^3/D_{4n}^*$ with
$kerf_*\cong\mathbb{Z}_{2n}$ and deg~$f=(2n)^2=4n^2\equiv 0
\text{\;(mod $4n$)}$.

When the self-mapping induces a endomorphism having the kernel as
$\mathbb{Z}_{\frac{2n}{(2n,\:k)}}$, according to the proof of lemma
3.1, $\frac{2n}{(2n,\:k)}$ should be an odd number. Take
$${\widetilde{f}}':S^3\to\,\mathbb{C}^2,(z_1,z_2)\mapsto
(z_1^{\frac{2n}{(2n,\:k)}},z_2^{\frac{2n}{(2n,\:k)}})$$ and let
$\widetilde{f}=\pi\,\circ{\widetilde{f}}':S^3\to\,S^3$. Then
$\widetilde{f}$ satisfy $\widetilde{f}\circ\,a=a\circ\widetilde{f}$,
$\widetilde{f}\circ\,b=\,b^{\frac{2n}{(2n,\:k)}}\circ\widetilde{f}$.
This induces $f:S^3/D_{4n}^*\to\,S^3/D_{4n}^*$ with
$kerf_*\cong\mathbb{Z}_{\frac{2n}{(2n,\:k)}}$ and
deg~$f$=\,$(\frac{2n}{(2n,\:k)})^2$ is any square of odd factors of
$n$.

Since the self-mapping-degrees can be multiplied by any
$k^2\;(\,(\,k\,,4n\,)=1\,)$ which is the degree of some self-mapping
inducing automorphism of $D_{4n}^*$\;(see \cite{Hayat4}),
\begin{displaymath}
  D(S^3/D_{4n}^*)=
  \{\,h^2\,|\,h\in\mathbb{Z};\,2\nmid\,h\;\text{or}\;h=n\;\text{or}\;h=0\}+4n\mathbb{Z}.
\end{displaymath}

\bigskip

\textbf{Case III:}

$\pi_1(M)=T_{24}^*=\langle\,a,b\,|\,a^2=b^3=(ab)^3=-1\,\rangle$\;.

Take $\rho:T_{24}^*\to\,SU(2)$ as $a\,\mapsto\left(
\begin{array}{ccc}
   i & 0\\
   0 & -i
\end{array} \right)$,

$b\,\mapsto\left(\begin{array}{ccc}
   \frac12+\frac12i & \frac12+\frac12i\\
  -\frac12+\frac12i & \frac12-\frac12i
\end{array}\right)=\left(\begin{array}{ccc}
     \frac{\sqrt2}2\zeta & \frac{\sqrt2}2\zeta\\
   \frac{\sqrt2}2\zeta^3 & \frac{\sqrt2}2\zeta^7
\end{array}\right)$,\,here\;$\zeta^4=-1$.

Then the left actions of elements of $T_{24}^*$ on $SU(2)$ are
$a\,\circ\,(z_1,z_2)=(iz_1,iz_2)$,
$b\,\circ\,(z_1,z_2)=(\frac{\sqrt2}2\zeta\,z_1-\frac{\sqrt2}2\zeta\,\bar{z}_2,
\frac{\sqrt2}2\zeta\,\bar{z}_1+\frac{\sqrt2}2\zeta\,z_2)$.

In lemma 3.1 we had proved that the non-trivial endomorphism of
$T_{24}^*$ has the kernel isomorphic to $Q_8$. Take
$$\widetilde{f}:S^3\to\,\mathbb{C}^2,(z_1,z_2)\mapsto
(z_1^4+\bar{z}_2^{\,4}+2iz_1^2\bar{z}_2^2,\;2\sqrt2\zeta\bar{z}_1^2z_2^2)$$
and let $\widetilde{f}=\pi\,\circ{\widetilde{f}}':S^3\to\,S^3$. Then
$\widetilde{f}\circ\,a=\widetilde{f},\widetilde{f}\circ\,b=b^2\circ\widetilde{f}$,
which induces $f:S^3/T_{24}^*\to\,S^3/T_{24}^*.$~
$Kerf_*\cong\,Q_8$.\;deg$f$=$16$.

In \cite{Hayat4} we know that the set of degrees of self-mappings
inducing isomorphisms of $T_{24}^*$ is $\{\,1\,\}+24\mathbb{Z}$.
Hence,
$$D(S^3/T_{24}^*)=\{\,0,1,16\}+24\mathbb{Z}.$$

\textbf{Case IV:}

$\pi_1(M)=O_{48}^*=\langle\,a,b\,|\,a^4=b^3=(ab)^2=-1\rangle$\;.

According to Lemma 2.2, there exists some mapping
$f:S^3/O_{48}^*\to\,S^3/O_{48}^*$ with $ker\,f_*=T_{24}^*$,
$f_*:a\mapsto-1,\;b\mapsto1$.

For the Sylow 2-subgroup
$H=\langle\,a,ba^2b^{-1}\,\rangle\cong\,D_{4\cdot\,4}^*$ of
$O_{48}^*$, construct covering map
$p:\widetilde{_HM}\cong\,S^3/D_{4\cdot\,4}^*\to\,S^3/O_{48}^*$ such
that $p_*\pi_1(\widetilde{_HM})=H$. Then $f$ can be lifted to
$\widetilde{_Hf}:\widetilde{_HM}\to\widetilde{_HM}$:
$$\CD
  \widetilde{_HM} @> \widetilde{_Hf} >> \widetilde{_HM} \\
  @V p VV @V p VV  \\
  S^3/O_{48}^* @>f>> S^3/O_{48}^*.
\endCD$$

Now $\pi_1(\widetilde{_HM})$ is $D_{4n}^*$ type and
$ker(\widetilde{_Hf})_*\cong\,Q_8\cong\,D_{4\cdot\,2}^*$
$=\langle\,a^2,\,ba^2b^{-1}\,\rangle$. We had discussed this case in
Case\;II. So ~deg~$\widetilde{_Hf}=4^2\equiv\,0\;\text{(mod\;$16$)}$
and ~deg~$f$~$\equiv\,0\;\text{(mod\;$16$)}$.

For the Sylow 3-subgroup $K=\langle\,b^2\,\rangle\cong\mathbb{Z}_3$,
construct the covering map $q:\widetilde{_KM}\to\,S^3/O_{48}^*$ and
lift $f$ to $\widetilde{_Kf}:\widetilde{_KM}\to\widetilde{_KM}$. It
is easy to check $(\widetilde{_Kf})_*$ is a constant endomorphism.
So ~deg~$f\equiv\,0\;\text{(mod\;3)}$. Together with
~deg~$f$~$\equiv\,0\;\text{(mod\;$16$)}$ we can get
~deg~$f\equiv\,0\;\text{(mod\;48)}$.

So non-trivial endomorphisms of $O_{48}^*$ give no new self-mapping
degree of $S^3/O_{48}^*$ other than the trivial ones. Hence,
together with the results of the degrees yielded by mappings
inducing automorphisms of $O_{48}^*$\,(\,see \cite{Hayat4}), we have

$$D(S^3/O_{48}^*)=\{\,0,1,25\}+48\mathbb{Z}.$$

\textbf{Case V:}

$\pi_1(M)=I_{120}^*=\langle\,a,b\,|\,a^2=b^3=(ab)^5=-1\rangle$\;.

$I_{120}^*$ has no non-trivial endomorphisms and the set of degrees
of self-mapping inducing isomorphisms is
$\{1,49\}+120\mathbb{Z}$\,(see\;\cite{Hayat4}). Hence,
$$D(S^3/I_{120}^*)=\{\,0,1,49\}+120\mathbb{Z}.$$

\textbf{Case VI:}

$\pi_1(M)=T_{8\cdot3^q}^{'}=\,Q_8\rtimes\mathbb{Z}_{3^q}$,\,
$Q_8=\{\,\pm\,1,\,\pm\,i,\,\pm\,j,\,\pm\,k\,\},\,\mathbb{Z}_{3^q}=\langle\,w\,\rangle$,
$w\cdot\,i\cdot\,w^{-1}=j\,;\;\;w\cdot\,j\cdot\,w^{-1}=k\,;\;\;w\cdot\,k\cdot\,w^{-1}=i.$

If a self-mapping $f$ of $S^3/T_{8\cdot\,3^q}'$ induces an
endomorphism of $T_{8\cdot\,3^q}'$ with the kernel
$K=Q_8\times\langle\,w^{3^{^{p}}}\rangle$, then
$f_*(Q_8)=1,\;f_*(w)=w^{\,k\cdot\,3^{^{q-p}}}\;(\,3\nmid\,k,
q\geq\,p>0\,)$. Its lifting
$(\widetilde{_{Q_8}f}):\widetilde{_{Q_8}M}\to\widetilde{_{Q_8}M}$\;($\widetilde{_{Q_8}M}$
is the covering space of $M$ with the fundamental group equals
$Q_8$) will satisfy $(\widetilde{_{Q_8}f})_*=1$ and its lifting
$\widetilde{_{Z_{3^q}}f}:\widetilde{_{Z_{3^q}}M}\to\widetilde{_{Z_{3^q}}M}\cong\,S^3/\mathbb{Z}_{3^q}$
will satisfy
$(\widetilde{_{Z_{3^q}}f})_*:w\mapsto\,w^{\,k\cdot\,3^{^{q-p}}}$. So
\[\left\{\begin{array}{ll}
  \text{deg}\,f\equiv\text{deg}\,\widetilde{_{Q_8}f}\equiv0        &(\text{mod\;}8)\\
  \text{deg}\,f\equiv\text{deg}\,\widetilde{_{Z_{3^q}}f}\equiv(k\cdot\,3^{q-p}\,)^2\qquad(\,3\nmid\,k\,)\;&(\text{mod}\;3^q).
\end{array}\right.\]
The solution is
\[
\text{deg}\,f\equiv\left\{\begin{array}{ll}
k^2\cdot\,(3^{2q-2p}-3^{q})&(\,2\,|\,q\,)\\
k^2\cdot\,(3^{2q-2p}-3^{q+1})&(\,2\nmid\,q\,)\\
\end{array}\right.
(\text{mod}\;8\cdot\,3^q)\quad(\,3\nmid\,k,\;q\geq\,p>0\,).
\]
After multiplied by a square of integer coprime to $8\cdot\,3^q$
which is a feasible degree of some mapping inducing an automorphism
of $\pi_1(M)$, the possible self-mapping degrees remain the same.
Hence,
$$
D(S^3/T_{8\cdot3^q}^{'})=\left\{\begin{array}{ll}
  \{\,k^2\cdot\,(3^{2q-2p}-3^{q})\;\;|\,3\nmid\,k, q\geq\,p>0\}+8\cdot3^q\mathbb{Z}&(\,2\,|\,q\,)\\
  \{\,k^2\cdot\,(3^{2q-2p}-3^{q+1})\;\;|\,3\nmid\,k,q\geq\,p>0\}+8\cdot3^q\mathbb{Z}&(\,2\nmid\,q\,).
\end{array}\right.
$$

\textbf{Case VII:}

$\pi_1(M)=D_{n'\cdot\,2^q}'=\,\mathbb{Z}_{n'}\rtimes\mathbb{Z}_{2^q}$,\,
$\mathbb{Z}_{n'}=\langle\,u\,\rangle,\,\mathbb{Z}_{2^q}=\langle\,w\,\rangle,\;w\cdot\,u\cdot\,w^{-1}=u^{-1}.$

Given $f:S^3/D_{n'\cdot\,2^q}'\to\,S^3/D_{n'\cdot\,2^q}'$ with
$f_*=\phi:u\mapsto1,\,w\mapsto\,w^{\,h\cdot\,2^{^{q-p}}}\;(\,2\nmid\,h,\,q\geq\,p>0\,)$,
the liftings
$\widetilde{_{\mathbb{Z}_{n'}}f}:\widetilde{_{\mathbb{Z}_{n'}}M}\to\widetilde{_{\mathbb{Z}_{n'}}M}$
and
$\widetilde{_{\mathbb{Z}_{2^q}}f}:\widetilde{_{\mathbb{Z}_{2^q}}M}\to\widetilde{_{\mathbb{Z}_{2^q}}M}$
will satisfy $(\widetilde{_{\mathbb{Z}_{n'}}f})_*:u\mapsto1$ and
$(\widetilde{_{\mathbb{Z}_{2^q}}f})_*:w\mapsto\,w^{\,h\cdot\,2^{^{q-p}}}$;
Given $f:S^3/D_{n'\cdot\,2^q}'\to\,S^3/D_{n'\cdot\,2^q}'$ with
$f_*=\phi:u\mapsto\,u^{l\cdot\,n'/n''},\,w\mapsto1\;(n''>1,
(l,n'')=1\,)$, these liftings will satisfy
$(\widetilde{_{\mathbb{Z}_{n'}}f})_*:u\mapsto\,u^{l\cdot\,n'/n''}$
and $(\widetilde{_{\mathbb{Z}_{2^q}}f})_*:w\mapsto1$.

The equations
\[
\left\{\begin{array}{ll}
  \text{deg}\,f\equiv\text{deg}\,\widetilde{_{\mathbb{Z}_{n'}}f}\equiv0      &(\text{mod\;}n')\\
  \text{deg}\,f\equiv\text{deg}\,\widetilde{_{\mathbb{Z}_{2^{^q}}}f}\equiv\,(h\cdot\,2^{^{q-p}})^{^2}
  \qquad(\,2\nmid\,h, q\geq\,p>0) &(\text{mod}\;2^{\,q})
\end{array}\right.\;
\]
imply
$$\text{deg}f\equiv\,h^2\cdot\,2^{\,2(q-p)}\cdot\,[\,1-2^{(\,2p-q)(n'-1)}\,]\;\;(\text{mod}\;n'\cdot\,2^{\,q}),$$

and the equations
\[
\left\{\begin{array}{ll}
  \text{deg}\,f\equiv\text{deg}\,\widetilde{_{\mathbb{Z}_{n'}}f}\equiv\,(l\cdot\,n'/n'')^2\qquad(\,(l,n'')=1)&(\text{mod\;}n'),\\
  \text{deg}\,f\equiv\text{deg}\,\widetilde{_{\mathbb{Z}_{2^{^q}}}f}\equiv0&(\text{mod}\;2^{\,q})
\end{array}\right.
\]
imply
$$\text{deg}f\equiv\,(l\cdot\,n'/n'')^2\,[\,1-(n')^{2^q-1}\,]\;\;(\text{mod}\;n'\cdot\,2^{\,q}).$$

After considering the composition of the above mappings, we know
that the self-mapping degree should be the multiplication of such
numbers. Hence
\[
  \begin{array}{ll}
    D(S^3/{D_{n'\cdot\,2^q}'})=&\{\,k^2\cdot[\,1-(n')^{2^{^{\,q}}-1}\,]^i\cdot
      ([\,1-2^{\,(\,2p-q\,)(\,n'-1\,)}\,])^j\;|\\
    &\qquad\qquad\,i,j,k,p\in\mathbb{Z},q\geq\,p>0\}+n'\,2^{\,q}\,\mathbb{Z}\,.
  \end{array}
\]

\bigskip

\textbf{Case VIII:}

$\pi_1(M)=\mathbb{Z}_m\times\,G$,\,$(m,|G|)=1$, $G$ is some group of
the above seven cases.

For a self mapping $f$ of $M$, construct the covering space of
$\widetilde{_GM}$ and $\widetilde{_{\mathbb{Z}_m}M}$ correspondence
to subgroups $G$ and $\mathbb{Z}_m$ of $\pi_1(M)$. Similarly,
consider the liftings
$\widetilde{_Gf}:\widetilde{_GM}\to\widetilde{_GM}$ and
$\widetilde{_{\mathbb{Z}_m}f}:\widetilde{_{\mathbb{Z}_m}M}\to\widetilde{_{\mathbb{Z}_m}M}$.
We have
\[\left\{\begin{array}{lll}
  \text{deg}\,f\equiv\text{deg}\,\widetilde{_Gf}\equiv\,a  & (a\in\,D(S^3/G))  &(\text{mod\;}|G|)\\
  \text{deg}\,f\equiv\text{deg}\,\widetilde{_{\mathbb{Z}_m}f}\equiv\,k^2 & (k\in\mathbb{Z})
  &(\text{mod}\;m).
\end{array}\right.\]
When $G$ is $D_{n'\cdot\,2^q}'$ and $T_{8\cdot3^q}^{'}$, the
solution do not have neat expressions. Nevertheless, the set of self
mapping degree of $M$ can be always completely determined by the
above equations.

\end{document}